\documentclass[preprint]{imsart}
\RequirePackage{amsthm,amsmath,amssymb,bbm}
\RequirePackage[]{natbib}
\RequirePackage[colorlinks,citecolor=blue,urlcolor=blue]{hyperref}

\startlocaldefs
\numberwithin{equation}{section}
\theoremstyle{plain}
\newtheorem{theorem}{Theorem}
\newtheorem{assume}{Assumption}

\newtheorem{cor}{Corollary}

\newtheorem{rmk}{Remark}

\newcommand{\cb}{{\cal B}}

\newcommand{\cK}{\mathcal{K}}

\newcommand{\argmin}{\mathop{\arg\min}}

\newcommand{\ud}{\mathrm{d}}

\newcommand{\bY}{{\bf Y}}
\newcommand{\bx}{{\bf x}}

\newcommand{\mE}{\mathbb{E}}

\newcommand{\be}{\begin{equation}}
	\newcommand{\ee}{\end{equation}}
\newcommand{\beqn}{\begin{eqnarray}}
	\newcommand{\eeqn}{\end{eqnarray}}
\endlocaldefs

\allowdisplaybreaks

\begin{document}
	\begin{frontmatter}
\title{Kullback–Leibler excess risk bounds for exponential weighted aggregation in Generalized linear models}
		\runtitle{Kullback–Leibler excess risk bounds for EWA}
\thankstext{T1}{The views, results, and opinions expressed in this work are solely those of the author and do not, in any way, represent those of the Norwegian Institute of Public Health.
		}
		
\begin{aug}
\author{\fnms{The Tien}~\snm{Mai}\ead[label=e1]{the.tien.mai@fhi.no}\orcid{0000-0002-3514-9636}}
			\address{
Norwegian Institute of Public Health, Oslo, Norway.
				\\
				\printead[presep={\ }]{e1}
			}
			\runauthor{T.T. Mai}
		\end{aug}
		
		\begin{abstract}
Aggregation methods have emerged as a powerful and flexible framework in statistical learning, providing unified solutions across diverse problems such as regression, classification, and density estimation. In the context of generalized linear models (GLMs), where responses follow exponential family distributions, aggregation offers an attractive alternative to classical parametric modeling. This paper investigates the problem of sparse aggregation in GLMs, aiming to approximate the true parameter vector by a sparse linear combination of predictors. We prove that an exponential weighted aggregation scheme yields a sharp oracle inequality for the Kullback-Leibler risk with leading constant equal to one, while also attaining the minimax-optimal rate of aggregation. These results are further enhanced by establishing high-probability bounds on the excess risk.
		\end{abstract}
		
		
		
	\end{frontmatter}

\section{Introduction}
\label{sc_introduction}
The problem of aggregation has attracted growing attention in the statistical learning community, as it provides a versatile framework capable of encompassing various learning scenarios. Initially introduced in the context of regression by \citet{nemirovski2000topics} and \citet{juditsky2000functional} as an extension of model selection, aggregation has since evolved into a mature area of research. The works of \citet{tsybakov2003optimal} and \citet{yang2004aggregating} were instrumental in this development, establishing optimal aggregation rates. Subsequent studies extended aggregation methods to density estimation \citep{rigollet2007linear} and classification \citep{belomestny2007spatial}, further illustrating their generality. Additional foundational contributions include \citet{catoni2004statistical}, \citet{leung2006information}, \citet{bunea2007aggregation}, and \citet{tsybakov2011exponential}, particularly in the context of Gaussian regression.

Generalized Linear Models (GLMs) offer a natural extension of classical Gaussian linear regression by allowing the response variable to follow a distribution from the exponential family, rather than being restricted to the normal distribution. This generalization accommodates a broad spectrum of practical applications, including binomial and Poisson data. The foundational framework for GLMs was laid out in the seminal work of \citet{mccullagh1989generalized}.

In this paper, we investigate the problem of aggregation in the GLM setting. More specifically, let $(\bx_i, Y_i)$, for $i = 1, \ldots, n$, denote observed data, where the distribution of $Y_i$ belongs to the exponential family with natural parameter $\theta^0_i$. Unlike standard GLMs, which assume a model of the form $\theta^0_i = \bx_i^\top \beta$, the aggregation framework avoids this modeling assumption. Instead, the objective is to approximate the true parameter vector $\theta^0 = (\theta^0_1, \ldots, \theta^0_n)$ by a linear combination $\theta_\beta = \sum_{j=1}^p \beta_j \bx_j$ that minimizes the Kullback-Leibler (KL) divergence, $KL(\theta^0, \theta_\beta)$, over a constraint set $\mathcal{B} \subseteq \mathbb{R}^p$.

The choice of $\mathcal{B}$ gives rise to different aggregation schemes, following the terminology of \citet{bunea2007aggregation} and \citet{abramovich2016model}:
\begin{itemize}
    \item {\em linear} aggregation
($\mathcal{B}=\mathcal{B}_L=\mathbb{R}^p$),

\item  {\em convex} aggregation
($\mathcal{B}=\mathcal{B}_C=\{\beta \in \mathbb{R}^p: \beta_j \geq 0,\;
\sum_{j=1}^p \beta_j=1\}$),

\item  {\em model selection} aggregation
($\mathcal{B}=\mathcal{B}_{MS}$ is a subset of vectors with a single nonzero entry),
and

\item  {\em subset selection} or {\em
	$p_0$-sparse} aggregation ($\mathcal{B}=\mathcal{B}_{SS}(p_0)=\{\beta \in
\mathbb{R}^p: ||\beta||_0 \leq p_0\}$ for a given $  p_0 \geq 1 $).
\end{itemize}
Notably, linear and model selection aggregation correspond to the extreme cases of subset selection, i.e., $\mathcal{B}_L = \mathcal{B}_{SS}(p) $ and $\mathcal{B}_{MS} = \mathcal{B}_{SS}(1)$. The problem of linear, convex, and model selection aggregation in GLMs was addressed in \citet{rigollet2012kullback} via maximum likelihood aggregation methods. The more general case of $p_0$-sparse aggregation in high dimensions was studied by \citet{abramovich2016model}, who proposed penalized maximum likelihood estimators.

Since the true parameter $\theta^0$ is unknown, the goal is to construct an estimator $\theta_{\widehat\beta}$ that closely mimics the ideal (oracle) approximation. Specifically, we aim to achieve a bound of the form \begin{equation} \label{eq:excess} \mathbb{E}[KL(\theta^0, \theta_{\widehat\beta})] \leq C \inf_{\beta \in \mathcal{B}} KL(\theta^0, \theta_\beta) + \Delta_{\mathcal{B}}(\theta^0, \theta_{\widehat\beta}), \end{equation} where $C \geq 1$ and $\Delta_{\mathcal{B}}$ is the excess KL risk. Ideally, one would like $C = 1$ and $\Delta_{\mathcal{B}}$ as small as possible. The paper \citet{rigollet2012kullback} established optimal asymptotic rates for $\Delta_{\mathcal{B}}$ and constructed estimators achieving these rates with $C = 1$ for linear, convex, and model selection aggregation. In contrast, while \citet{abramovich2016model} showed that their penalized likelihood estimator achieves minimax-optimal rates for subset selection aggregation, they obtained only a non-sharp oracle inequality with $C \geq 4/3$.

When the model is misspecified—i.e., $\inf_{\beta \in \mathcal{B}} KL(\theta^0, \theta_\beta) > 0$—obtaining $C = 1$ becomes crucial. \citet{abramovich2016model} conjectured that achieving such sharp oracle inequalities would require model averaging, rather than selecting a single model. In this work, we confirm this conjecture. By employing an exponential weighting aggregation (EWA) procedure, we derive a sharp oracle inequality with $C = 1$ and provide bounds on the excess KL risk that hold with high probability.

EWA has attracted considerable interest in both statistics and machine learning communities, see for example \citep{guedj2019primer,alquier2024user,hellstrom2025generalization}. It has been successfully applied to various problems such as sparse regression \citep{dalalyan2008aggregation, alquier2011PAC, dalalyan2012sparse}, classification \citep{catonibook, mai2023high}, matrix problems \citep{dalalyan2020exponential, mai2023reduced}, and deep learning \citep{tinsi2022risk, mai2025misclassification}. In this paper, we adapt the EWA methodology to the GLM setting.

Specifically, we derive a sharp oracle inequality of the form \eqref{eq:excess} with leading constant $C = 1$, while also attaining the optimal aggregation rate established in \cite{abramovich2016model}. In addition, we extend their results by providing high-probability bounds on the excess risk. To exploit the underlying sparsity structure, we adopt a scaled Student’s t-distribution prior on the coefficients—a prior that has been effectively utilized in various sparse estimation problems \cite{dalalyan2012mirror, dalalyan2012sparse, mai2023high, mai2024sparse}, but has not yet been applied within the GLM framework.

The remainder of the paper is structured as follows. In Section \ref{sc_model_method}, we introduce the problem setup and present our sparse EWA methodology. Section \ref{sc_theory} contains our main theoretical results. All technical proofs are deferred to Section \ref{sc_proof}.

\section{Setup and notation}
\label{sc_model_method}
\subsection{Setup}
Consider a GLM setup with a response variable $Y$ and a vector of $ p $ dimensional predictor $x $. We observe a series of independent
observations $(\bx_i,Y_i),\;i=1,\ldots, n$, where the design points $\bx_i \in
\mathbb{R}^p$ are deterministic. The distribution $f_{\theta_i}(y)$ of $Y_i$
belongs to a one-parameter natural exponential family with a natural
parameter $\theta_i$ and a scaling parameter $a$:
\begin{equation}
  \label{eq:model} 
f_{\theta_i}(y)=\exp\left\{\frac{y \theta_i -
	b(\theta_i)}{a}+c(y,a)\right\}
    .  
\end{equation}
Here, $a$ is assumed to be known.
The function $b(\cdot)$ is assumed to be twice-differentiable.
In this case, 
$\mE(Y_i)=b'(\theta_i)$ and $Var(Y_i)=ab''(\theta_i)$ (see \cite{mccullagh1989generalized}).

To fully specify a generalized linear model (GLM), we adopt the canonical link function, which expresses the natural parameter for each observation as $\theta_i = \bx_i^t\beta $. In matrix notation, this relationship is written as $\theta_\beta = X\beta$, where $X \in \mathbb{R}^{n \times p}$ is the design matrix and $\beta \in \mathbb{R}^p$ represents the vector of unknown regression coefficients. In this paper, we assume that the design matrix is fixed.

\subsection{Kullback--Leibler excess risk}
\label{subGAP}

The Kullback--Leibler divergence between two
probability distributions $P$ and $Q$ is defined by
\[
\cK(P\|Q)= \begin{cases}
\int\log\bigl(\frac{\ud P}{\ud Q} \bigr)\,\ud P, &\quad {\rm if }\,
	P \ll Q,\vspace*{2pt}\cr
	\infty, &\quad {\rm otherwise}.
\end{cases}
\]

Let $f_{\theta}$ and $f_{\zeta}$ be two possible joint distributions of
the data from the exponential family with $n$-dimensional vectors of natural parameters
$\theta$ and $\zeta$ correspondingly.
A Kullback-Leibler divergence $ \mathcal{K} (\theta \| \zeta) $ between
$f_{\theta}$ and $f_{\zeta}$ is then
\be
\mathcal{K} (\theta \| \zeta)
    =
    \mE_{\theta}\left\{ \ln
	\left(\frac{f_{\theta}(\bY)}{f_{\zeta}(\bY)}\right)\right\}
    =
    \frac{1}{a}~\mE_{\theta}
	\left\{\sum_{i=1}^n Y_i(\theta_i-\zeta_i)-b(\theta_i)+b(\zeta_i)\right\} 
\ee
where $b(\theta)=(b(\theta_1),\cdots,b(\theta_n))$ and $b(\zeta)=(b(\zeta_1),\ldots,b(\zeta_n))$.

For a given estimator $\widehat{\theta}$ of the unknown $\theta$ consider
$  \mathcal{K} (\theta \| \widehat{\theta}) $, the Kullback-Leibler divergence between the
true distribution $f_{\theta}$ of the data and 
its empirical distribution $f_{\widehat{\theta}}$ generated by $\widehat{\theta}$. The goodness of $\widehat{\theta}$ is measured
by the corresponding Kullback-Leibler \textit{risk}:
\be \label{eq:mE}
\mE \, [\mathcal{K} (\theta \| \widehat{\theta})]
=
\frac{1}{a}\left[ b'(\theta)^t(\theta-\mE(\widehat{\theta}))-(b(\theta)-\mE b(\widehat{\theta}))^t
{\bf 1}\right]
\ee
where the expectation is taken w.r.t. the true distribution $f_{\theta}$.
In particular, for the Gaussian case,
where $b(\theta)=\theta^2/2$ and $a=\sigma^2$, $\mE [\mathcal{K} (\theta \| \widehat{\theta})] $ is the mean squared error
$ \mE ||\widehat{\theta}-\theta||^2$ divided by the constant $2\sigma^2$.

The quantity of interest is the following  Kullback-Leibler \textit{excess} risk:
%
%
\begin{equation}
	\label{EQKLag}
\mE \, [\mathcal{K} ( \theta^0 \| \theta_{\widehat{\beta} } )] 
-
\min_{\beta \in \mathcal{B}_{SS}(p_0) } 
 [\mathcal{K} (\theta^0 \| \theta_\beta )]
    ,
\end{equation}
where
$  \mathcal{B}_{SS}(p_0)
=
\{\beta \in
\mathbb{R}^p: ||\beta||_0 \leq p_0\}
$ for a given $ 1 \leq p_0 \leq p $. Here, $ \theta^0 $ is the true unknown parameter.
When working with large datasets, it is often assumed that only a limited number of predictors significantly influence the response. This  sparsity assumption is particularly important in high-dimensional settings where the number of predictors exceeds the number of observations ($p > n$) \cite{abramovich2016model, abramovich2018high}.

\subsection{A sparse EWA approach}

Putting
\begin{equation}
    r(\beta) 
    =
\frac{1}{n}  \frac{1}{a} \sum_{i=1}^n [ Y_i (X \beta)_i - b( (X \beta)_i) ]
    .
\end{equation}
Note that the above empirical risk function \( r(\beta) \) proportional to the log-likelihood. Moreover, in this context, the maximum likelihood estimator coincides with the minimizer of the Kullback-Leibler risk defined in \eqref{eq:mE}; see also \cite{rigollet2012kullback, abramovich2016model}. Instead of pursuing a maximum likelihood-type aggregation, we adopt an exponentially weighted aggregation approach.

For any \(\lambda > 0\), with a sparsity-inducing prior distribution $ \pi $ given in \eqref{eq_priordsitrbution}, we consider the following exponentially weighted aggregation, \(\hat{\rho}_{\lambda} \),
\begin{equation}
\label{eq_Gibbs_poste}
\hat{\rho}_{\lambda} (\beta) 
\propto 
\exp[-\lambda r(\beta)] \pi(\beta)
,
\end{equation}

The goal of the above aggregation is to shift the distribution in favor of parameter values that achieve lower empirical risk on the observed data. The tuning parameter $\lambda$ controls the strength of this adjustment, and its role will be examined in more detail in the following sections. Specifically, if $ \pi $ assigns higher probability to sparse vectors, $\hat{\rho}_\lambda $ will favor sparse vectors with low empirical risk, satisfying our requirements.

The EWA shown in \eqref{eq_Gibbs_poste} is also often called the Gibbs posterior as in \cite{alquier2016,catonibook,dalalyan2012sparse,dalalyan2008aggregation}. The selection of $\hat{\rho}_\lambda $ is based on Donsker and Varadhan's variational formula, presented in Lemma \ref{lemma:dv}, rather than following traditional Bayesian approaches. Throughout this paper, we use $\pi$ to represent the prior and $\hat{\rho}_\lambda $ as the pseudo-posterior. Various applications of the EWA method are covered in recent reviews in \cite{guedj2019primer, alquier2024user}.

In this study, we examine the following prior distribution, proposed in \cite{dalalyan2012mirror,dalalyan2012sparse}.
For a fixed constant \( B_1 >0 \), for all $ \beta \in \mathbb{R}^p  $ that \(  \|\beta\|_2 \leq B_1 \), we adopt the following scaled Student distribution as our prior distribution,
\begin{eqnarray}
\label{eq_priordsitrbution}
\pi (\beta) 
\propto 
\prod_{i=1}^{p} 
(\zeta^2 + \beta_{i}^2)^{-2}
,
\end{eqnarray}
where $ \zeta>0 $ is a tuning parameter.   The constant \( B_1 \) is conventionally taken to be large, leading to an approximate distribution of \( \pi \) as \( T\zeta \sqrt{2} \), where \( T \) is a random vector whose components are independently drawn from a Student’s t-distribution with 3 degrees of freedom. Setting \( \zeta \) to a sufficiently small value ensures that most entries of \( \zeta T \) are concentrated near zero, while the heavy-tailed behavior of the distribution permits occasional large deviations. This structure promotes sparsity in the parameter vector by leveraging the prior distribution.

\section{Main results}
\label{sc_theory}
\subsection{Assumption}
We assume the following assumption on the parameter space $\Theta$
and the second derivative $b''(\cdot)$.

\begin{assume} 
\label{asssume_A}
 Assume that $\theta_i \in \Theta$, where the parameter space $\Theta \subseteq \mathbb{R}$ is a closed (finite or infinite) interval.
 Assume that there exist a constant $0 < U < \infty$ such that the function $b''(\cdot)$ satisfies that:
$\sup_{t \in \Theta} b''(t) \leq U $.
\end{assume}
Comparable assumptions have been made in the context of generalized linear models, as seen in \cite{van2008high,rigollet2012kullback,abramovich2016model}. The condition on \( b''(\cdot) \) in Assumption \ref{asssume_A} is specifically designed to rule out extreme scenarios in which the variance \( \text{Var}(Y) \) becomes unbounded.  For Gaussian distribution, $b''(\theta)=1$ and, therefore,
$ U = 1 $ for any $\Theta$. For the binomial distribution,
$b''(\theta)=\frac{e^\theta}{(1+e^\theta)^2},$ and thus $ U = 1/4 $.

Put
\begin{align}
\label{eq_true_beta_star}
    \beta^* = \argmin_{\beta \in \mathcal{B}_{SS}(p_0) } 
 \mathcal{K} (\theta^0 \| \theta_\beta )	
    .
\end{align}

\subsection{Bounds in expectation}

We begin by providing an upper bound on the excess risk associated with \eqref{EQKLag}. A more refined non-asymptotic bound that holds with high probability will be established later in Theorem \ref{thm_main}.

\begin{theorem}
\label{thm_expectation}
Assume that Assumption \ref{asssume_A} is satisfied. Let $ \lambda = n $ and $ \zeta = 1/(np\| X \|) $.  Then for $ \beta^* $ such that $  \| \beta^*\|_2 \leq B_1 - 2d\zeta $ we have:
	\begin{align}
\mathbb{E}~\mathbb{E}_{\beta \sim \hat{\rho}_\lambda}
  \mathcal{K}(\theta^0 \| \theta_\beta )  
    - 	
\min_{\beta \in \mathcal{B}_{SS}(p_0) } 
 \mathcal{K} (\theta^0 \| \theta_\beta )			
\leq 
\mathfrak{C} p_0 \log (\frac{ np\|X\| }{ p_0})  
	,
	\end{align}
where $ \mathfrak{C} > 0 $ is a universal constant depending only on $ U, a, B_1 $.
\end{theorem}

\begin{rmk}
In Theorem \ref{thm_expectation}, we show that the integrated Kullback–Leibler (KL) risk of our method closely matches the best achievable KL risk for subset selection aggregation over the class \( \mathcal{B}_{SS}(p_0) \). This is an important result, particularly in scenarios where the true model lies outside the sparse subset class—i.e., when \( \min_{\beta \in \mathcal{B}_{SS}(p_0)} \mathcal{K}(\theta^0 \| \theta_\beta) > 0 \). Even in such misspecified settings, our method remains competitive, as it achieves a KL risk comparable to that of the best possible sparse approximation within \( \mathcal{B}_{SS}(p_0) \). This highlights both the robustness and practical relevance of our approach.
\end{rmk}

\begin{rmk}
The result in Theorem \ref{thm_expectation} importantly answers an open question raised in \cite{abramovich2016model} as they only obtain 
	\begin{align*}
\mathbb{E}
  \mathcal{K} (\theta^0 \| \theta_{\widehat{\beta} } )  
    - 	
C \min_{\beta \in \mathcal{B}_{SS}(p_0) } 
 \mathcal{K} (\theta^0 \| \theta_\beta )			
\leq 
C' p_0 \log (\frac{ np}{ p_0})  
	,
	\end{align*}
for some constant $ C \geq 4/3 $ while we obtain similar results with constant $ C =1 $.
\end{rmk}

\begin{rmk}
Interestingly, Theorem \ref{thm_expectation} also resolves another open question posed in \cite{alquier2020concentration}. Specifically, that paper presents a non-sharp oracle inequality where the left-hand side involves the $\alpha$-Rényi divergence, and the right-hand side involves the Kullback-Leibler divergence. The authors were unable to derive a result similar to Theorem \ref{thm_expectation} using their approach, and they left the task of establishing an oracle inequality with the Kullback-Leibler divergence on the left-hand side as an open problem. Further details can be found in Theorem 2.7 and Section 6 of \cite{alquier2020concentration}.
\end{rmk}

The primary technical tool employed in our proofs is the PAC-Bayesian bound method, which provides a powerful framework for deriving non-asymptotic risk bounds. However, we adopt a specialized formulation that incorporates the Kullback–Leibler risk in combination with a concentration inequalitie tailored to the setting of GLMs, as introduced in \cite{rigollet2012kullback}. Our methodology is more closely aligned with the influential approach developed in \cite{catoni2003PAC,catoni2004statistical,catonibook}, which emphasizes oracle-type PAC-Bayesian inequalities. These results offer refined risk bounds that adapt to the complexity of the model class under consideration. It is worth noting that PAC-Bayesian bounds were initially introduced in \cite{mcallester1998some,STW} to analyze the generalization error of Bayesian-type estimators. Since then, this framework has been extensively studied and expanded upon in a variety of contexts. See \cite{guedj2019primer} and \cite{alquier2024user,hellstrom2025generalization} for recent reviews. 
In particular, PAC-Bayesian techniques have been used to establish oracle inequalities in several high-dimensional estimation problems. Some applications include sparse regression \cite{alquier2013sparse,alquier2011PAC} and low-rank modeling \cite{mai2015,mai2017pseudo,mai2023bilinear,mai2023reduced}. These developments highlight the flexibility and strength of PAC-Bayesian theory as a unifying tool for statistical learning.

An important special case of Theorem \ref{thm_expectation} arises when \( \min_{\beta \in \mathcal{B}_{SS}(p_0)} \mathcal{K}(\theta^0 \| \theta_\beta) = 0 \), indicating that there exists some \( \beta^0 \in \mathcal{B}_{SS}(p_0) \) such that \( \theta^0 = X\beta^0 \), i.e., the true underlying model is exactly sparse. In this setting, we immediately obtain the following corollary.

\begin{cor}
\label{cor_expect_truemodel}
Assume that Theorem \ref{thm_expectation} is satisfied and there exist a $ \beta^0 \in \mathcal{B}_{SS}(p_0) $ such that $ \theta^0 = X\beta^0 $, then we have:
	\begin{align}
\mathbb{E}~
  \mathcal{K}(\theta^0 \| \theta_{\hat{\beta}_M}  )  		
\leq 
\mathfrak{C} p_0 \log ( np\|X\| /  p_0)  
	,
	\end{align}
where $ \mathfrak{C} > 0 $ is a universal constant depending only on $ U, a, B_1 $.
\end{cor}

\begin{rmk}
We highlight that the excess Kullback–Leibler risk established in Theorem \ref{thm_expectation} is essentially of order \( p_0 \log (p / p_0) \). Notably, this rate is achieved without requiring prior knowledge of the true sparsity level \( p_0 \), which demonstrates the adaptive nature of our proposed method. 
Even more importantly, up to a log-term, this rate aligns with the known minimax lower bound for the KL risk in our setting, as established in Theorem 2 of \cite{abramovich2016model}. Thus, our result not only adapts to unknown sparsity but also achieves the optimal rate from an information-theoretic perspective.
\end{rmk}

\begin{assume} 
\label{asssume_B}
 Assume that $\theta_i \in \Theta$, where the parameter space $\Theta \subseteq \mathbb{R}$ is a closed (finite or infinite) interval.
 Assume that there exist a constant $0 < L \leq U < \infty$ such that the function $b''(\cdot)$ satisfies that:
$\inf_{t \in \Theta} b''(t) \geq L $
.
\end{assume}

\begin{theorem}[Theorem 2 in \cite{abramovich2016model}] 
\label{th:lower}
Consider a GLM with the canonical link $\theta=X\beta$
under Assumption \ref{asssume_A} and \ref{asssume_B}. Let $1 \leq p_0 \leq r$.
Then, there exists a constant $C_2>0$ such that
\be
\label{eq:lower}
\inf_{\widehat{\theta}} \sup_{\beta \in \cb(p_0)}
\mE \,   \mathcal{K} 
 ( \theta_\beta \| \widehat{\theta}) 
\geq
C_2~\frac{ L }{U }
 p_0 \log (\frac{pe}{p_0})
,
\end{equation}
where the infimum is taken over all estimators $\hat{\theta}$ of $\theta$.
\end{theorem}

Note that Assumption \ref{asssume_B} was also employed in \cite{van2008high,rigollet2012kullback,abramovich2016model}. It excludes the degenerate case where \( \text{Var}(Y) = 0 \) and further ensures the strong convexity of \( b(\cdot) \) over \( \Theta \). 

It is noted that the reference \cite{abramovich2016model} also requires this assumption to derive upper bounds on the excess Kullback–Leibler (KL) risk, we do not rely on it for our main results. Nevertheless, if we impose Assumption \ref{asssume_B}, we can immediately obtain the following result for the mean estimator, as the KL risk becomes convex in its second argument under this condition.

We define our mean estimator as
\begin{equation}
    \hat{\beta}_M := \int \beta  \hat{\rho}_{\lambda} (\beta) {\rm d} \beta
\end{equation}
We immediately obtain the following corollary for our mean estimator from Theorem \ref{thm_expectation}.

\begin{cor}
\label{cor_expect}
Assume that Theorem \ref{thm_expectation} is satisfied and in addition Assumption \ref{asssume_B} is hold, then we have:
	\begin{align}
\mathbb{E}~
  \mathcal{K}(\theta^0 \| \theta_{\hat{\beta}_M}  )  
    - 	
\min_{\beta \in \mathcal{B}_{SS}(p_0) } 
 \mathcal{K} (\theta^0 \| \theta_\beta )			
\leq 
\mathfrak{C} p_0 \log (\frac{ np\|X\| }{ p_0})  
	,
	\end{align}
where $ \mathfrak{C} > 0 $ is a universal constant depending only on $ U, a, B_1 $.
\end{cor}

\subsection{Bounds in high probability}
While Theorem \ref{thm_expectation} already represents a significant improvement over the results established in \cite{abramovich2016model}, it is possible to strengthen these findings further. In particular, we extend our analysis by deriving non-asymptotic bounds that hold with high probability, thereby offering stronger probabilistic guarantees. This refinement is formalized in the following theorem.

\begin{theorem}
\label{thm_main}
Assume that Assumption \ref{asssume_A} is satisfied. Let $ \lambda = n $ and $ \zeta = 1/(np \| X \| ) $.  Then for $ \beta^* $ such that $  \| \beta^*\|_2 \leq B_1 - 2d\zeta $, we have with probability at least $ 1-\varepsilon, \varepsilon\in (0,1) $:
	\begin{align}
\mathbb{E}_{\beta \sim \hat{\rho}_\lambda}
  	\mathcal{K}(\theta^0 \| \theta_\beta )  - 	
\min_{\beta \in \mathcal{B}_{SS}(p_0) } 
 \mathcal{K} (\theta^0 \| \theta_\beta )		
\leq 
\mathfrak{A} \bigg[ p_0 \log (\frac{ np\|X\| }{ p_0})
	+ \log(1/\varepsilon) \bigg]
	,
	\end{align}
where $ \mathfrak{A} > 0 $ is a universal constant depending only on $ U, a, B_1 $.
\end{theorem}

Theorem \ref{thm_main} establishes a non-asymptotic oracle inequality for our proposed method. Comparable bounds for the maximum likelihood aggregation approach in the contexts of linear, convex, and model selection aggregation were previously derived in \cite{rigollet2012kullback}. To the best of our knowledge, our result is the first to provide such a bound for the case of $p_0$-sparse aggregation.

Analogous to Corollary \ref{cor_expect}, and under the additional Assumption \ref{asssume_B}, the following corollary for our mean estimator directly follows from Theorem \ref{thm_main}.

\begin{cor}
Suppose that Theorem \ref{thm_main} is satisfied and, in addition, Assumption \ref{asssume_B} holds. Then, we have with probability at least $ 1-\varepsilon, \varepsilon\in (0,1) $ that
	\begin{align}
\mathcal{K}(\theta^0 \| \theta_{ \hat{\beta}_M } )  - 	
\min_{\beta \in \mathcal{B}_{SS}(p_0) } 
 \mathcal{K} (\theta^0 \| \theta_\beta )		
\leq 
\mathfrak{A} \bigg[ p_0 \log (\frac{ np\|X\| }{ p_0})
	+ \log(1/\varepsilon) \bigg]
	,
	\end{align}
where $ \mathfrak{A} > 0 $ is a universal constant depending only on $ U, a, B_1 $.
\end{cor}

When the true underlying model is exactly sparse—i.e., there exists \( \beta^0 \in \mathcal{B}_{SS}(p_0) \) such that \( \theta^0 = X\beta^0 \). Analogously to Corollary \ref{cor_expect_truemodel}, we immediately obtain the following corollary.

\begin{cor}
\label{cor_proba_truemodel}
Assume that Theorem \ref{thm_main} is satisfied and there exist a $ \beta^0 \in \mathcal{B}_{SS}(p_0) $ such that $ \theta^0 = X\beta^0 $. Then we have with probability at least $ 1-\varepsilon, \varepsilon\in (0,1) $:
	\begin{align}
\mathbb{E}_{\beta \sim \hat{\rho}_\lambda}
  	\mathcal{K}(\theta^0 \| \theta_\beta ) 	
\leq 
\mathfrak{A} \bigg[ p_0 \log (\frac{ np\|X\| }{ p_0})
	+ \log(1/\varepsilon) \bigg]
	,
	\end{align}
where $ \mathfrak{A} > 0 $ is a universal constant depending only on $ U, a, B_1 $.
\end{cor}

\begin{rmk}
Our proposed approach demonstrates robustness with respect to the selection of various sparsity-inducing priors. Specifically, it is not restricted to a single form of prior specification, which enhances its flexibility and broad applicability in different modeling scenarios. For instance, one can incorporate a model selection-type prior, such as the one introduced in \cite{alquier2011PAC}, which is well-suited for enforcing exact sparsity through a discrete prior structure. Alternatively, other continuous shrinkage priors like the Horseshoe prior \cite{carvalho2010horseshoe} may also be employed. In such cases, theoretical guarantees can still be obtained by leveraging recent results as in \cite{mai2024concentration}. This adaptability ensures that our method remains effective across a wide range of prior formulations commonly used in sparse Bayesian modeling.
\end{rmk}

It is noted that frequentist properties of Bayesian methods in the context of sparse GLMs have been studied in the literature. For example, \cite{jeong2021posterior} investigates the posterior contraction rates under sparsity assumptions, providing valuable insights into the asymptotic behavior of the posterior distribution. However, the focus of their work is fundamentally different from ours. While they concentrate on the rate at which the posterior contracts around the true parameter value in a frequentist sense, our primary interest lies in the analysis of Kullback–Leibler risk. This divergence in focus reflects distinct inferential goals: theirs rooted in estimation accuracy, and ours in predictive performance and information-theoretic guarantees.

\subsection{Example}

Here, we present an example of our results for the case of Gaussian distribution with known variance $ \sigma^2 $. In our model \eqref{eq:model}, we have that $a =\sigma^2 $ and $ b(\theta) = \theta^2/2 $ and our empirical risk function is as
\begin{equation*}
    r(\beta) 
    =
\frac{1}{n \sigma^2 } \sum_{i=1}^n [ Y_i (X \beta)_i -  (X \beta)_i^2/2 ]
    .
\end{equation*}
We immediately obtain the following results from Theorem \ref{thm_expectation} and \ref{thm_main} for the case of Gaussian distribution as follows. Note that in this case Assumption \ref{asssume_A} is always satisfied with constant $ U = 1$. Note also that in this case also have $  \mathcal{K}(\theta^0 \| \theta_\beta ) = \frac{1}{2\sigma^2} \| \theta^0  - \theta_\beta \|^2 $.

\begin{cor}
\label{thm_gaussian}
Let $ \lambda = n $ and $ \zeta = 1/(np \| X \| ) $.  Then for $ \beta^* $ such that $  \| \beta^*\|_2 \leq B_1 - 2d\zeta $, we have with probability at least $ 1-\varepsilon, \varepsilon\in (0,1) $:
	\begin{align}
\mathbb{E}_{\beta \sim \hat{\rho}_\lambda}
\| \theta^0  - \theta_\beta \|^2 
- 	
\!\!
\min_{\beta \in \mathcal{B}_{SS}(p_0) } 
\| \theta^0 
 - \theta_\beta \|^2	
\leq 
\mathfrak{G} \bigg[ p_0 \log (\frac{ np\|X\| }{ p_0})
+ \log(\frac{1}{\varepsilon}) \bigg]
,
\end{align}
where $ \mathfrak{G} > 0 $ is some universal constant depending only on $ \sigma, B_1 $.
\end{cor}

Similar bounds in expectation for sparse Gaussian regression have been established for a related estimator in \cite{dalalyan2008aggregation}. In contrast, our work not only provides high-probability bounds but also extends these results to the broader class of generalized linear models.

\section{Proof}
\label{sc_proof}

\subsection*{Acknowledgments}
The results, views, and opinions presented in this paper are solely those of the author and do not, in any way, represent those of the Norwegian Institute of Public Health.

\subsubsection*{Conflicts of interest/Competing interests}
The author declares no potential conflict of interests.

\end{document}